\documentclass[11pt,fullpage]{article}

\RequirePackage{amsthm,amsmath,amsfonts,amssymb}
\RequirePackage[numbers]{natbib}

\RequirePackage[colorlinks,citecolor=blue,urlcolor=blue]{hyperref}
\RequirePackage{graphicx}

\theoremstyle{plain}

\newtheorem{theorem}{Theorem}[section]

\theoremstyle{definition}

\newtheorem{assumption}[theorem]{Assumption}

\theoremstyle{remark}

\numberwithin{equation}{section}

\long\def\ignore#1{}

\newcommand{\be}{\begin{equation}}
\newcommand{\ee}{\end{equation}}
\newcommand{\beqn}{\begin{eqnarray}}
\newcommand{\eeqn}{\end{eqnarray}}

\newcommand{\bfm}[1]{\mbox{\boldmath{$#1$}}}
\newcommand{\balpha}{\bfm{\alpha}}
\newcommand{\bbeta}{\bfm{\beta}}

\newcommand{\blambda}{\bfm{\lambda}}
\newcommand{\bkappa}{\bfm{\kappa}}
\newcommand{\bPsi}{\bfm{\Psi}}

\newcommand{\bx}{{\bf x}}
\newcommand{\bX}{{\bf X}}

\newcommand{\bm}{{\bf m}}

\newcommand{\mE}{\mathbb{E}}
\newcommand{\mP}{\mathbb{P}}
\newcommand{\cA}{{\cal A}}
\newcommand{\cB}{{\cal B}}
\newcommand{\cG}{{\cal G}}

\newcommand{\cJ}{{\cal J}}

\newcommand{\cm}{{\cal M}}

\newcommand{\cE}{{\cal E}}
\newcommand{\cF}{{\cal F}}
\newcommand{\cP}{{\cal P}}
\newcommand{\cC}{{\cal C}}

\newcommand{\cX}{{\cal X}}

\renewcommand{\top}{T}

\begin{document}

\newcommand{\bb}[1]{{\mathbb #1}}

\long\def\ignore#1{}

\title{\bf Classification by sparse generalized additive models}
\author{{\bf Felix Abramovich} \\
Department of Statistics and Operations Research \\
Tel Aviv University \\
Israel\\
felix@tauex.tau.ac.il
}

\date{}

\maketitle

\begin{abstract}
We consider (nonparametric) sparse (generalized) additive models (SpAM) for classification.
The design of a SpAM classifier is based on minimizing the logistic loss with a sparse group Lasso/Slope-type penalties on the coefficients of univariate additive components' expansions in orthonormal series (e.g., Fourier or wavelets). The resulting classifier is inherently adaptive to the unknown sparsity and smoothness. We show that under certain sparse group restricted eigenvalue condition it is  nearly-minimax (up to log-factors) simultaneously across the entire range of analytic, Sobolev and Besov classes.
The performance of the proposed classifier is illustrated on a simulated and a real-data examples. 
\end{abstract}

\section{Introduction} \label{sec:introduction}
Consider a binary classification setup with a vector of features $\bX \in \mathbb{R}^d$ and the outcome class label $Y|(\bX=\bx) \sim B(1,p(\bx))$, where $p(\bx)=P(Y=1|\bX=\bx)$. 
The quality of a classifier $\eta: \mathbb{R}^d \rightarrow \{0,1\}$ is quantified by its misclassification error $R(\eta)=P(Y \neq \eta(\bX))$ and the optimal Bayes classifier is
$\eta^*(\bx)=I\{p(\bx) \geq 1/2\}$. However, the conditional probability $p(\bx)$ is typically unknown. One then estimates it from the data $D=(\bX_1,Y_1),\ldots,(\bX_n,Y_n)$ by some $\widehat{p}(\bx)$ and designs a plug-in classifier $\widehat{\eta}=I\{\widehat{p}(\bx) \geq 1/2\}$ with the conditional risk $R(\widehat{\eta})=P(Y \neq \widehat{\eta}(\bX)|D)$. The goodness of $\widehat{\eta}$ is measured by the misclassification excess risk $\cE(\widehat{\eta},\eta^*)=\mE R(\widehat{\eta})-R(\eta^*)$, where the expectation is w.r.t. the distribution of $D$.

Probably the most commonly used model for $p(\bx)$ is logistic regression, where it is assumed that the log-odds (logit)
\be \label{eq:logistic}
g(\bx)=\ln \frac{p(\bx)}{1-p(\bx)}=\bbeta^\top \bx
\ee
and $\bbeta \in \mathbb{R}^d$ is a vector of (unknown) regression coefficients. The corresponding Bayes (linear) classifier for the model (\ref{eq:logistic}) is $\eta^*(\bx)=I\{\bbeta^\top \bx \geq 0\}$.
Assuming $d < n$ and estimating $\bbeta$ by the maximum likelihood estimator (MLE) $\widehat{\bbeta}$ yields the 
plug-in classifier $\widehat{\eta}_L(\bx)=I\{\widehat{\bbeta}^\top \bx \geq 0\}$. Logistic regression classification has been  well-studied and is one of the main tools used by practitioners. In particular, it is well-known known that
$\cE(\widehat{\eta}_L,\eta^*) \sim \sqrt{d/n}$ which  is rate-optimal (in the minimax sense) over all linear classifiers (see, e.g. \cite{ag19}).

Yet, a linear model for $g(\bx)$ in (\ref{eq:logistic}) might be too restrictive and oversimplified in a variety of applications. Nonparametric models allow one more flexibility by assuming much more general assumptions on $g(\bx)$. One can use then various nonparametric techniques (e.g., kernel estimators, local polynomial regression, orthogonal series or nearest neighbours) to estimate the unknown $p(\bx)$ (or, equivalently, $g(\bx)$) and define a nonparametric plug-in classifier $\widehat{\eta}_N(\bx)=I\{\widehat{p}(\bx) \geq 1/2\}$. \cite{y99} showed that with the properly chosen tuning parameters, $\widehat{\eta}_N(\bx)$ achieves the minimax rates across various smoothness classes.
Thus, assuming that the unknown $p(\bx)$ has a smoothness $s$, $\cE(\widehat{\eta}_N,\eta^*)=O(n^{-\frac{s}{2s+d}})$.

A general nonparametric classifier $\widehat{\eta}_N$ suffers, however, from a ``curse of dimensionality'' problem, where the convergence rates rapidly slow down as $d$ increases and the required sample size for consistent classification grows exponentially with $d$. 
To handle it one needs some additional assumptions on the model.

Note that classification is particularly challenging near the decision boundary $\{\bx: p(\bx)=1/2\}$ or, equivalenlty, $\{\bx: g(\bx)=0\}$, where it is especially hard to predict the class label accurately, while at points far from the boundary, it is easier. Thus, assume that for all $0<h<h^*$, $P(|p(\bX)-1/2| \leq h) \leq C h^\alpha$ for some $C>0, \alpha \geq 0$ and $0<h^*<1/2$ known as 
the margin or low-noise condition \cite{t04}. The case $\alpha=0$ essentially corresponds to no assumption, while the second extreme case $\alpha=\infty$ implies the existence of a hard margin of size $h^*$ separating between two classes. \cite{at07} showed that under the above margin assumption, the rates of misclassification excess risk of $\widehat{\eta}_N$ can be indeed improved. Assuming again that $p(\bx)$ has a smoothness $s$, under some mild conditions on the marginal distribution $\mP_X$ of $\bX$, the rate for the misclassification excess risk becomes of the order $n^{-\frac{s(1+\alpha)}{(2+\alpha)s+d}}$ approaching the parametric rate $n^{-1}$ as $\alpha \rightarrow \infty$.
It can even be reduced further under stronger assumptions on $\mP_X$.

Although the margin assumption improves the convergence rates, it provides only a partial solution to 
the curse of dimensionality problem  since computational burden for fitting a general nonparametric classifier remains prohibitive even for relatively small and moderate $d$.
In addition, it is difficult to interpret such a model, to perform feature selection, etc.

To overcome these challenges one needs additional structural assumptions on the model. A common approach is to consider (generalized) {\em additive} models (GAM)(\cite{ht99}). 
In particular, for logistic regression, assume that
\be \label{eq:additive}
g(\bx)=\mu+\sum_{j=1}^d g_j(x_j),
\ee
where $g_j(\cdot)$'s are univariate ``smooth'' functions. To make the model (\ref{eq:additive}) identifiable  impose $\mE g_j(X_j)=0$ for all $j=1,\ldots,d$. The GAM (\ref{eq:additive}) is a natural nonparametric generalization of the linear logistic regression model (\ref{eq:logistic}). Additive models have become a standard tool in high-dimensional nonparametric regression and classification and can be efficiently fitted by the backfitting algorithm \cite{ht99}.
Additivity assumption  drastically improves the convergence rates of the resulting classifier $\widehat{\eta}_A$. Thus, adapting the results of \cite{y99} and \cite{rwy12} one
can show that $\cE(\widehat{\eta}_A,\eta^*)=O(\sqrt{d}~ n^{-\frac{s}{2s+1}})$ if all univariate components $g_j$ are of smoothness $s$, or, more generally, $\cE(\widehat{\eta}_A,\eta^*)=O\left(\sqrt{\sum_{j=1}^d n^{-\frac{2s_j}{2s_j+1}}}\right)$ if $g_j$ are of different smoothness.

Nevertheless, in the era of ``Big Data'', the number of features $d$ might be very large and even larger that the sample size $n$
(large $d$, small $n$ setups). The GAM (\ref{eq:additive})  requires $d \ll n$ and cannot cope with a curse of dimensionality in this case. 
Reducing dimensionality of a feature space by selecting a sparse subset of significant active features becomes essential. Thus, consider a {\em sparse} generalized additive model (SpAM) assuming that
\be \label{eq:SpAM}
g(\bx)=\mu+\sum_{j \in \cJ} g_j(x_j),
\ee
where $\cJ \subseteq \{1,\ldots,d\}$ is an (unknown) subset of active features of cardinality $|\cJ|=d_0$ that have a ``significant'' impact on the outcome. We are mainly interested in sparse setups, where $d_0 \ll \min(d,n)$.

SpAMs have been intensively studied in the context of nonparametric regression w.r.t. the quadratic risk, where a common approach is based on penalized least squares estimation  with different combinations of sparsity-induced and regularized smoothness convex penalties.
See, e.g., \cite{lz06, mgb09, rllw09, ky10, rwy12, al15}  for various SpAM regression estimators and their properties. Their comparison is discussed in \cite{al15} and \cite{hss22}. \cite{ky10} and \cite{hss22} extended this approach to more general convex losses relevant, in particular, for the considered logistic regression model. However, their estimators are not adaptive to the (usually unknown) smoothness of $g_j$'s.

%We consider SpAM within classification framework.
Classification by sparse linear logistic regression was studied in \cite{ag19}. They showed that the corresponding minimax misclassification excess risk is of the order $\sqrt{d_0 \ln(de/d_0)/n}$ and designed penalized MLE classifiers that achieve it.

In this paper we investigate (nonparametric) SpAM classifiers. We propose
a SpAM estimator of $g$ in (\ref{eq:SpAM}) defined in the Fourier/wavelet domain w.r.t. the penalized logistic loss with {\em sparse group Lasso} or more general {\em sparse group Slope} penalties on the Fourier/wavelet coefficients. 
The estimator and the resulting plug-in classifier are inherently {\em adaptive} to the unknown sparsity and smoothness classes. 
We establish the minimax rates for the misclassification excess risk for SpAMs across analytic, Sobolev and Besov classes of functions and show that there is a phase transition between sparse and dense SpAMs. 
In particular, if $g_j,\;j \in \cJ$ has a smoothness $s_j$, the minimax misclassificitaion excess risk is of the order $\max\left(\sqrt{\frac{d_0 \ln(de/d_0)}{n}}, \sqrt{\sum_{j \in \cJ} n^{-\frac{2s_j}{2s_j+1}}}\right)$. 
We prove that with a proper choice of tuning parameters and under certain restricted sparse group eigenvalue condition,
the sparse group Lasso/Slope SpAM classifiers $\widehat{\eta}_{sgL}$ and $\widehat{\eta}_{sgS}$
are simultaneously nearly-minimax (up to log-factors) across the entire range of those classes.

The rest of the paper is organized as follows. Section \ref{sec:model} presents the SpAM model for classification and defines sparse group Lasso and Slope classifiers. In Section \ref{sec:minimaxity} we establish the minimaxity of the proposed SpAM classifiers across Sobolev, analytic and Besov classes.
Their performance is illustrated on a simulated and a real-data examples in Section \ref{sec:example}. All the proofs are given in Appendix \ref{sec:appendix}.

\section{SpAM classifier} \label{sec:model}
Consider a sparse additive logistic regression model:
\be \label{eq:model}
Y|(\bX=\bx) \sim B(1,p(\bx)),
\ee
where $\bX \in \mathbb{R}^d$ is a vector of linearly independent features with a marginal probability distribution $\mP_X$ on a bounded support $\cX$. Without loss of generality assume that $\cX \subseteq [0,1]^d$. 
%Let $V=\mE(\bX \bX^\top)$ be the second moment matrix of $\bX$.
The logit function $g(\bx)=\ln\frac{p(\bx)}{1-p(\bx)} \in L_2(\cX)$ has a sparse additive form (\ref{eq:SpAM}), where $g_j$ lies in some class of functions $\cF_j$ on the unit interval (e.g., Sobolev $H^{s_j}[0,1]$ or more general Besov $B^{s_j}_{p_j,q_j}[0,1]$), $\mE g_j(X_j)=0$ and the subset of active features $\cJ$ is a priori unknown.

Expand the unknown univariate functions $g_j,\;j=1,\ldots,d$ in some orthonormal basis 
$\{\psi_\ell\}_{l=0}^\infty$ in $\cF_j$ as
$$
g_j(x_j)=\sum_{\ell=0}^\infty \beta_{j \ell} \psi_\ell(x_j),
$$
where $\beta_{j\ell}=\int g_j(x_j)\psi_\ell(x_j)dx_j$. For simplicity of exposition we consider the same basis for all $g_j$'s. Due to the identifiability conditions $\mE g_j(X_j)=0$, we consider $\beta_{j0}=0$ for all $j$. By Parseval's identity
$||g_j||^2_{L_2[0,1]}=\sum_{l=1}^\infty |\beta_{j\ell}|^2$ and $||g||^2_{L_2[0,1]^d}=\mu^2+\sum_{j=1}^d \sum_{l=1}^\infty |\beta_{j\ell}|^2$. 

Estimate the unknown coefficients $\beta_{j\ell}$ from the data sample $(\bX_i,Y_i),\;i=1,\ldots,n$.
Obviously, for all $j=1,\ldots,d$, we set $\widehat{\beta}_{j\ell}=0$ for all $\ell > n$. To estimate $\beta_{j\ell}$ for $1 \leq \ell \leq n$ consider additive truncated estimators of the form 
$$
\widetilde{g}_n(\bx)=\widetilde{\mu}+\sum_{j=1}^d \sum_{\ell=1}^n \widetilde{\beta}_{j\ell}\psi_\ell(x_j)
$$
and the logistic loss 
\be \label{eq:logloss}
\begin{split}
\ell(Y,\widetilde{g}_n(\bX)) & =\ln\left(1+\exp(\widetilde{g}_n(\bX))\right) -Y \widetilde{g}_n(\bX) \\
&=\ln \left(1+\exp\left(\widetilde{\mu}+tr(\Psi(\bX)^\top \widetilde{B})\right)\right)-\left(\widetilde{\mu}+tr(\Psi(\bX)^\top \widetilde{B})\right) Y,
\end{split}
\ee
where $\widetilde{B} \in \mathbb{R}^{d \times n}$ is the regression coefficients matrix and $\Psi(\bX) \in \mathbb{R}^{d \times n}$ is the  matrix with entries $\Psi_{j\ell}(\bX)=\psi_\ell(X_j)$. 

Let $B \in \mathbb{R}^{d \times n}$ be the matrix with the true (unknown) regression coefficients $\beta_{j\ell},\;j=1,\ldots,d;\; \ell=1,\ldots n$.
Sparse representation of $g$ and smoothness properties of $g_j$'s in (\ref{eq:SpAM}) are directly expressed in terms of $B$.
Evidently, all $\beta_{j\ell}=0$ for all $\ell$ for $j \notin \cJ$. In addition, if the chosen basis $\{\psi_\ell\}_{l=0}^\infty$ allows sparse representation of
$g_j, \; j \in \cJ$, such that they can be well-approximated by a small number of basis functions, then the decreasingly ordered $|\beta|_{j(\ell)}$ tend rapidly to zero as $\ell$ increases for  $j \in \cJ$. These arguments induce a double row-wise sparsity structure of $B$: it has only $d_0$ non-zero rows (global row-wise sparsity) and even those are ``approximately sparse'' in the sense that they have only a few ``large'' entries (local row-wise approximate sparsity). To capture such type of sparsity we consider  a {\em logistic sparse group Lasso} estimator of $B$. Sparse group Lasso was proposed in \cite{sfht13} for linear regression. Its logistic version was used in multiclass classification by multinomial linear logistic regression in \cite{vh14} and \cite{la23}. Within the considered SpAM binary classification setup we define the logistic sparse group Lasso estimator $\widehat{B}_{sgL}$ of $B$ as follows:
\be \label{eq:sparsegroupLasso}
\begin{split}
(\widehat{B}_{sgL},\widehat{\mu}_{sgL})=\arg \min_{\widetilde{B},\widetilde{\mu}} & \left\{\frac{1}{n}\ell(Y_i,\widetilde{g}_n(\bX_i))+\lambda \sum_{j=1}^d |\widetilde{B}_{j\cdot}|_2+\kappa \sum_{j=1}^d |\widetilde{B}_{j \cdot}|_1  \right\} \\
= \arg \min_{\widetilde{B},\widetilde{\mu}} & \left\{\frac{1}{n}\sum_{i=1}^n \left(\ln \left(1+\exp\left(\widetilde{\mu}+tr(\Psi(\bX_i)^\top \widetilde{B})\right)\right) \right. \right.\\ 
&-\left. \left. \left(\widetilde{\mu}+tr(\Psi(\bX_i)^\top \widetilde{B})\right) Y_i\right)+ 
\lambda \sum_{j=1}^d |\widetilde{B}_{j\cdot}|_2+\kappa \sum_{j=1}^d |\widetilde{B}_{j \cdot}|_1 \right\},
\end{split}
\ee
where $|\widetilde{B}_{j\cdot}|_2$ and $|\widetilde{B}_{j\cdot}|_1$ are respectively the $l_2$ and $l_1$-norms of a $j$-th row of $\widetilde{B}$, and $\lambda>0$ and $\kappa>0$ are tuning parameters.

The corresponding estimator 
\be \label{eq:g_sgL}
\widehat{g}_{sgL}(\bx)=\widehat{\mu}_{sgL}+\sum_{j=1}^d \sum_{l=1}^n \widehat{\beta}_{sgL} \psi_\ell(x_j)
\ee
and the resulting logistic sparse group Lasso classifier 
\be \label{eq:SpAMlasso}
\widehat{\eta}_{sgL}(\bx)=I\{\widehat{g}_{sgL}(\bx) \geq 0\}.
\ee

A more flexible generalization of (\ref{eq:sparsegroupLasso}) is a 
{\em logistic sparse group Slope} estimator that utilizes {\em sorted} norms in the penalty:
\be \label{eq:sparsegroupSlope}
\begin{split}
\widehat{B}_{sgS}= \arg \min_{\widetilde{B},\widetilde{\mu}} 
& \left\{ \frac{1}{n}\sum_{i=1}^n\left(\ln \left(1+ \exp\left(\widetilde{\mu}+tr(\Psi(\bX_i)^\top \widetilde{B})\right)\right)\right. \right. \\
& - \left. \left. \left(\widetilde{\mu}+tr(\Psi(\bX_i)^\top \widetilde{B})\right) Y_i\right) \right.  \\
& + \left. 
\sum_{j=1}^d \lambda_j |B|_{(j)2}+\sum_{j=1}^d \sum_{\ell=1}^n \kappa_\ell |\widetilde{B}|_{j(\ell)}  \right\},
\end{split}
\ee
where with some ambiguity of notations $|\widetilde{B}|_{(1)2} \geq \ldots \geq |\widetilde{B}|_{(d)2}$ are the descendingly ordered $l_2$-norms of the rows of $\widetilde{B}$, $|\widetilde{B}|_{j(1)} \geq \ldots \geq |\widetilde{B}|_{j(n)}$ are the descendingly ordered absolute values of entries of its $j$-row, and $\lambda_1 \geq \ldots \geq \lambda_d>0$ and $\kappa_1 \geq \ldots \geq \kappa_n>0$ are sequences of tuning parameters.  
Evidently, $\widehat{B}_{sgL}$ is a particular case of $\widehat{B}_{sgS}$
with equal $\lambda_j$'s and $\kappa_\ell$'s.
\ignore{
A similar form of logistic sparse group Slope classifier was proposed for multiclass multinomial linear logistic regression classification in \cite{la23}.
}

Similarly to sparse group Lasso classification, we consider 
$$
\widehat{g}_{sgS}(\bx)=\widehat{\mu}_{sgS}+\sum_{j=1}^d \sum_{l=1}^n \widehat{\beta}_{sgS} \psi_\ell(x_j)
$$
and the corresponding sparse group Slope classifier
\be \label{eq:SpAMslope}
\widehat{\eta}_{sgS}(\bx)=I\{\widehat{g}_{sgS}(\bx) \geq 0\}.
\ee

\section{Minimaxity across various function classes} \label{sec:minimaxity}
We now show that with the proper choice of tuning parameters the proposed logistic sparse group Lasso and Slope classifiers are adaptively nearly-minimax (up to log-factors) across various classes of functions $\cF_j$ for $g_j$'s in (\ref{eq:SpAM}).

As usual with convex penalization, one needs some (mild) assumptions on the design. Consider a truncated version $g_n$ of $g$:
$$
g_n(\bx)=\mu+\sum_{j=1}^d \sum_{\ell=1}^n \beta_{j\ell} \psi_\ell(x_j).
$$
Let random vectors $\bPsi(X_j)=(\psi_1(X_j),\ldots,\psi_n(X_j))^\top,\; j=1,\ldots,d$ be the rows of the matrix $\Psi(\bX)$ and define the  cross-covariance matrices 
$$
V_{jk}=E_X (\bPsi(X_j)\bPsi(X_k)^\top). 
$$
By simple algebra,
$$
\int g_n(\bx)^2d\mP_X=\mu^2+\sum_{j=1}^d\sum_{k=1}^d B_{j\cdot}V_{jk}B_{k\cdot}^\top.
$$

For a given $1 \leq d_0 \leq d$, denote a vector of integers $\bm=(m_1,\ldots,m_{d_0}) \in \mathbb{R}^{d_0}$ with $m_j \in \{1,\ldots,n\}$.
\begin{assumption}{($WRE(d_0,\bm,c_0)$-condition).} \label{as:design}
Given non-increasing sequences $\{\lambda_j\}_{j=1}^d$ and $\{\kappa_\ell\}_{\ell=1}^n$, consider the cone of matrices 
$$
{\cal S}(d_0,\bm,c_0)=\left\{A \in \mathbb{R}^{d \times n}:
||A||_{\blambda,\bkappa} \leq c_0
\left(\sqrt{\sum_{j=1}^{d_0} \lambda_j^2}+ \sqrt{\sum_{j=1}^{d_0} \sum_{\ell=1}^{m_j} \kappa_\ell^2}\right)||A||_F\right\},
$$
where $||A||_{\blambda,\bkappa}=\sum_{j=1}^d \lambda_j |A|_{(j)2}+\sum_{j=1}^{d} \sum_{\ell=1}^n \kappa_\ell |A|_{j(\ell)}$ is the sparse group Slope norm of $A$ and $||A||_F$ is its Frobenius norm,
and assume that
$$
    \nu_{gS}(d_0)= \inf_{A \in {\cal S}(d_0,\bm,c_0)\backslash \{0_{n \times d}\}}
    \frac{\sum_{j=1}^d \sum_{k=1}^d A_{j\cdot}V_{jk}A_{k\cdot}^\top}{||A||^2_F} > 0 .
    $$
\end{assumption}
The weighted restricted eigenvalue (WRE) condition \ref{as:design} is an extension of the WRE condition introduced for Slope in \cite{blt18}  to sparse group Slope and random design. Such type of assumption is common for convex penalties
(see \citep{blt18} for discussion).
In particular, for isotropic $\bX$, due to orthogonality of the basis, $V_{jj}=I_n$ and  $V_{jk}=0_{n \times n}$ for $j \neq k$. Thus,
$\sum_{j=1}^d \sum_{k=1}^d A_{j\cdot}V_{jk}A_{k\cdot}^\top=
||A||^2_F$ and Assumption
\ref{as:design} is trivially satisfied. 

Evidently, for sparse group Lasso with constant $\lambda$ and $\kappa$
\begin{equation} \nonumber
\begin{split}
{\cal S}(d_0,\bm,c_0)=\left\{A \in \mathbb{R}^{d \times n}: \right. & \lambda \sum_{j=1}^d |A_{j\cdot}|_2+\kappa\sum_{j=1}^n |A_{j\cdot}|_1 
\leq   \\
& \left. c_0
\left(\lambda \sqrt{d_0}+\kappa \sqrt{\sum_{j=1}^{d_0} m_j}\right) ||A||_F\right\}.
\end{split}
\end{equation}

\subsection{Sobolev classes} \label{subsec:sobolev}
Consider the orthonormal cosine series $\psi_0(x_j)=1,\;\psi_{\ell}(x_j)=\sqrt{2}\cos(\pi \ell x_j),\;\ell=1,2,\ldots$ and let $\beta_{j\ell}=\sqrt{2}\int_0^1 g_j(t)\cos(\pi \ell t)dt$
be the cosine Fourier coefficients of $g_j$. 

Suppose that $g_j$ belongs to a periodic Sobolev  class \newline $\widetilde{H}^{s_j}[0,1]=\left\{g_j: \sum_{\ell=0}^\infty \ell^{2s_j} \beta^2_{j\ell} \leq R_j,\;s_j>\frac{1}{2}\right\}$. In particular, for an integer $s_j$ it is equivalent to
$$\left\{g_j: \int_0^1 \left(g_j^{(s_j)}(t)\right)^2dt \leq R_j,\;g_j^{(m)}(0)=g_j^{(m)}(1)=0,\;m=1,\ldots,s_j-1\right\}$$ 
(e.g., \cite{t09}, Section 1.7.1).

Consider a set of sparse additive functions
\begin{equation} \nonumber
\begin{split}
\cG_{\widetilde H}(d_0,{\bf s})= & \left\{g(\bx) \in L_2(\cX):~g(\bx)=\mu+\sum_{j \in \cJ} g_j(x_j),~|\cJ| \leq d_0, \right. \\
& \left. g_j \in \widetilde{H}^{s_j}[0,1],\; j \in \cJ,\;\mE g_j(X_j)=0 \right\}
\end{split}
\end{equation}
and the corresponding set of logistic SpAM classifiers
$$
\cC_{\widetilde H}(d_0,{\bf s})=\left\{\eta(\bx)=I\{g(\bx) \geq 0\}: g(\bx) \in \cG_{\widetilde H}(d_0,{\bf s})\right\}.
$$

The following theorem establishes the upper bounds for the misclassification excess risk for logistic sparse group Lasso and Slope classifiers over periodic Sobolev classes:
\begin{theorem} \label{th:sobolevupper}
Consider the SpAM  (\ref{eq:logistic})-(\ref{eq:SpAM}), where $g \in \cG_{\widetilde H}(d_0,{\bf s})$, and assume $WRE(d_0,\bm,c_0)$-condition \ref{as:design} with $m_j=n^{\frac{1}{2s_j+1}},\;j=1,\ldots,d_0$ and $c_0$ that can be derived from the proof.

\begin{enumerate}
\item Apply the sparse group Lasso classifier (\ref{eq:SpAMlasso}) with $\lambda=C_1 \sqrt{\frac{\ln d}{n}}$
and $\kappa=C_2 \sqrt{\frac{\ln n}{n}}$ for some $C_1, C_2>0$ that can be derived from the proof.
Then, 
$$
\sup_{\eta^* \in \cC_{\widetilde H}(d_0, {\bf s})}\cE(\widehat{\eta}_{sgL},\eta^*) 
\le C \sqrt{\frac{1}{\nu_{sg}(d_0)}~\left(\frac{d_0 \ln d}{n}+
\ln n \sum_{j \in \cJ} n^{-\frac{2s_j}{2s_j+1}}\right)}
$$
for some $C>0$.

\item Apply the sparse group Slope classifier (\ref{eq:SpAMslope}) with $\lambda_j=C_1 \sqrt{\frac{\ln(de/j)}{n}},\;j=1,\ldots,d$
and $\kappa_\ell=C_2 \sqrt{\frac{\ln(ne/\ell)}{n}},\;\ell=1,\ldots,n$. 
Then, 
$$
\sup_{\eta^* \in \cC_{\widetilde H}(d_0, {\bf s})}\cE(\widehat{\eta}_{sgS},\eta^*) 
\le C \sqrt{\frac{1}{\nu_{sg}(d_0)}~\left(\frac{d_0 \ln\left(\frac{de}{d_0}\right)}{n}+
\ln n \sum_{j \in \cJ} n^{-\frac{2s_j}{2s_j+1}}\right)}.
$$
\end{enumerate}
\end{theorem}
The proposed SpAM classifiers are inherently adaptive to a sparsity parameter $d_0$, a set of active features $\cJ$ and smoothness parameters $s_j$'s.
Theorem \ref{th:sobolevupper} implies that a more flexible Slope-type penalty allows one to reduce the log-factor in the first term of the upper bounds. 
We now show that the upper bounds in Theorem \ref{th:sobolevupper} are nearly-minimax (up to log-factors) over $\cC_{\widetilde H}(d_0, {\bf s})$:
\begin{theorem} \label{th:sobolevlower}
Consider the SpAM  (\ref{eq:logistic})-(\ref{eq:SpAM}), where $d_0\ln(\frac{de}{d_0}) \leq n$ and $g \in \cG_{\widetilde H}(d_0,{\bf s})$. Then,
$$
\inf_{\widetilde{\eta}}\sup_{\eta^* \in \cC_{\widetilde H}(d_0, {\bf s}),\mP_X}\cE(\widetilde{\eta},\eta^*) \geq \tilde{C} \sqrt{\frac{d_0 \ln\left(\frac{de}{d_0}\right)}{n}+
\sum_{j \in \cJ} n^{-\frac{2s_j}{2s_j+1}}}
$$
for some $\tilde{C}>0$, where the infimum is taken over all classifiers $\widetilde{\eta}$ based on the data $(\bX_1,Y_1),\ldots,(\bX_n,Y_n)$.
\end{theorem}

The obtained misclassification excess risk bounds contain two terms and indicate on a phase transition between sparse and dense SpAMs. The term $\frac{d_0 \ln(de/d_0)}{n} \sim \frac{1}{n}\ln \binom{d}{d_0}$ corresponds to the error of selecting $d_0$ nonzero univariate components $g_j$'s out of $d$ and is common for model selection. The term
$\sum_{j \in \cJ} n^{-\frac{2s_j}{2s_j+1}}$ is due to simultaneous nonparametric estimating $d_0$ functions $g_j \in \widetilde{H}^{s_j}([0,1])$.
Assume for simplicity that $s_j=s$ for all $j \in \cJ$, so that the second term becomes $d_0  n^{-\frac{2s}{2s+1}}$. One then immediately realizes that the model selection error is dominating for sparse cases, where $\frac{d_0}{d} \lesssim e^{-n^{\frac{1}{2s+1}}}$, while for less sparse cases nonparametric estimation is the main error source. A similar phase transition phenomenon for sparse  linear classification was shown in \cite{ag19}. 

It can be shown that the exact minimax risk (without the extra $\ln n$-factor) can be achieved by the penalized SpAM estimator of \cite{hss22} applied to the classification setup. However, it is not adaptive to smoothness, as both penalties and optimal values of the smoothing parameters depend on $s_j$'s. 
For a particular case of design on a regular lattice,
\cite{al15} proposed an adaptive SpAM estimator based on a certain complexity penalty on the number of nonzero rows of $B$ and their nonzero entries, with the exact minimax risk within the Gaussian additive models framework.  Sparse group Lasso/Slope can be viewed as a convex surrogate of their complexity penalty. We conjecture that under certain conditions the results of \cite{al15} can be extended to a general design and the logistic loss, but the use of a complexity penalty will make it computationally infeasible in this case.

\subsection{Analytic functions} \label{subsec:analytic}
Assume now that $g_j$'s are very smooth, where $g_j$ belongs to a class of analytic functions $\cA^{\alpha_j}[0,1]=\left\{g_j:  \sum_{\ell=0}^\infty e^{2\alpha_j \ell} \beta^2_{j\ell} \leq R_j,\;\alpha_j>0\right\}$. Such functions admit analytic continuation into a band of width $\alpha_j$ of the complex plane (e.g., \cite{glt96}) and, in particular, are infinitely differentiable.

Define a set of functions
\begin{equation} \nonumber
\begin{split}
\cG_{\cA}(d_0,\balpha)= & \left\{g(\bx) \in L_2(\cX):~g(\bx)=\mu+\sum_{j \in \cJ} g_j(x_j),~ |\cJ| \leq d_0, \right.\\
& \left. g_j \in \cA^{\alpha_j}[0,1],\; j \in \cJ,\;\mE g_j(X_j)=0\right\}
\end{split}
\end{equation}
and the corresponding set of SpAM classifiers
$$
\cC_\cA(d_0,\balpha)=\left\{\eta(\bx)=I\{g(\bx) \geq 0\}: g(\bx) \in \cG_{\cA}(d_0,\balpha)\right\}.
$$

We now show that the logistic sparse group Lasso classifier $\widehat{\eta}_{sgL}$ in (\ref{eq:SpAMlasso}) applied to the Fourier cosine coefficients $\beta_{j\ell}$'s with the same choice for tuning parameters as in Theorem \ref{th:sobolevupper} for Sobolev classes is also nearly-minimax 
over $\cC_\cA(d_0,\balpha)$,
while the logistic sparse group Slope classifier $\widehat{\eta}_{sgS}$ in (\ref{eq:SpAMslope}) 
achieves the exact minimax rate for these classes.

\begin{theorem} \label{th:analyticupper}
Consider the SpAM  (\ref{eq:logistic})-(\ref{eq:SpAM}), where $ g \in \cG_{\cA}(d_0,\balpha)$, and assume $WRE(d_0,\bm,c_0)$-condition \ref{as:design} with $m_j=\frac{1}{2\alpha_j}\ln\left(\frac{n}{\ln n}\right),\;j=1,\ldots,d_0$ and $c_0$ that can be derived from the proof. 

\begin{enumerate}
\item 
Apply the sparse group Lasso classifier (\ref{eq:SpAMlasso}) with $\lambda=C_1 \sqrt{\frac{\ln d}{n}}$
and $\kappa=C_2 \sqrt{\frac{\ln n}{n}}$ for some $C_1, C_2>0$ that can be derived from the proof.
Then, 
$$
\sup_{\eta^* \in \cC_{\cA}(d_0, \balpha)}\cE(\widehat{\eta}_{sgL},\eta^*) 
\le C \sqrt{\frac{1}{\nu_{sg}(d_0)}~\left(\frac{d_0 \ln d}{n}+
\sum_{j \in \cJ}\frac{1}{\alpha_j}~ \frac{\ln n}{n}\right)}
$$
for some $C>0$.

\item Apply the sparse group Slope classifier (\ref{eq:SpAMslope}) with $\lambda_j=C_1 \sqrt{\frac{\ln(de/j)}{n}},\;j=1,\ldots,d$
and $\kappa_\ell=C_2 \sqrt{\frac{\ln(ne/\ell)}{n}},\;\ell=1,\ldots,n$. 
Then, 
$$
\sup_{\eta^* \in \cC_{\cA}(d_0, \balpha)}\cE(\widehat{\eta}_{sgS},\eta^*) 
\le C \sqrt{\frac{1}{\nu_{sg}(d_0)}~\left(\frac{d_0 \ln\left(\frac{de}{d_0}\right)}{n}+
\sum_{j \in \cJ}\frac{1}{\alpha_j}~ \frac{\ln n}{n}\right)}.
$$
\end{enumerate}
\end{theorem}

\begin{theorem} \label{th:analyticlower}
Consider the SpAM  (\ref{eq:logistic})-(\ref{eq:SpAM}), where $d_0\ln(\frac{de}{d_0}) \leq n$ and $g \in \cA(d_0,\balpha)$. Then,
$$
\inf_{\widetilde{\eta}}\sup_{\eta^* \in \cC_{\widetilde H}(d_0, \balpha),\mP_X}\cE(\widetilde{\eta},\eta^*) \geq \tilde{C} \sqrt{\frac{d_0 \ln\left(\frac{de}{d_0}\right)}{n}+
\sum_{j \in \cJ}\frac{1}{\alpha_j}~ \frac{\ln n}{n}}
$$
for some $\tilde{C}>0$.
\end{theorem}
Similar to Sobolev classes there is the phase transition phenomenon for the misclassification excess risk. Assuming for simplicity equal $\alpha_j=\alpha$, the model selection error term $\frac{d_0 \ln(de/d_0)}{n}$ dominates
for sparse cases where $\frac{d_0}{d} \lesssim  n^{-\frac{1}{\alpha}}$.

\subsection{Besov classes} \label{subsec:besov}
Consider now more general Besov classes $B^s_{p,q}$ of functions. The precise formal definition of Besov spaces can be found, e.g. in \cite{m92}. On the intuitive level, (not necessarily integer) $s$ measures the number of function's derivatives, where their existence is required in the $L_p$-sense, while $q$ provides a further finer gradation. The Besov spaces include, in particular,
the traditional Sobolev $H^{s}$ and H\"older $C^s$ classes of smooth functions
($B^s_{2,2}$ and $B^{s}_{\infty,\infty}$ respectively) but also various classes of spatially inhomogeneous functions like functions of bounded variation, sandwiched between $B^1_{1,1}$ and $B^1_{1,\infty}$ (\cite{m92, dj95}).

Besov classes can be equivalently characterized by coefficients of function's expansion in orthonormal wavelet series. Given a compactly supported scaling function $\phi$ of regularity $r$ and the corresponding mother wavelet $\psi$, 
one generates an orthonormal wavelet basis on the unit interval with $\psi_{hk}(t)=2^{h/2}\psi(2^h t-k),\;h \geq 0,\;k=0,\ldots,2^h-1$ (see \cite{cdv93}).  
Assume that $g_j \in B^{s_j}_{p_j,q_j}[0,1],\;1 \leq p_j,q_j \leq \infty,\; (\frac{1}{p_j}-\frac{1}{2})_+ < s_j < r$ and let $\beta_{j,hk}=\int_0^1 g_j(t) \psi_{hk}(t)dt$ be its wavelet coefficients. Then, 
\be \label{eq:besov}
\sum_{h=0}^\infty \left(2^{h(s_j+1/2-1/p_j)}\left(\sum_{k=0}^{2^h-1} |\beta_{j,hj}|^{p_j}\right)^{1/{p_j}}\right)^{q_j} \leq R_j
\ee
with $l_{p_j}$ and/or $l_{q_j}$
norms in (\ref{eq:besov}) being replaced by the corresponding $l_\infty$-norms for $p_j=\infty$ and/or $q_j=\infty$ 
\ignore{
For example
$$
\sup_{0 \leq h \leq \infty} 2^{h+1/2-1/p_j}\left(\sum_{k=0}^{2^h-1} |\beta_{j,hj}|^{p_j}\right)^{1/{p_j}} \leq R_j,\;\;\;1 \leq p_j <\infty, q_j=\infty
$$
}
(e.g., \cite{m92, dj95}). 

Moreover, let $\ell=2^h+k,\;h \geq 0,\;0 \leq k \leq 2^h-1$ and re-write the triple-indexed wavelet coefficients $\beta_{j,hk}$ in terms of double-indexed $\beta_{j\ell}$. Then,
\be \label{eq:weaklp}
|\beta|_{j(\ell)} \leq C \ell^{-(s_j+1/2)}
\ee
(\cite{d93}, Lemma 2).

Define the set of sparse additive functions
\begin{equation} \nonumber
\begin{split}
\cG_{\cB}(d_0,{\bf s}, {\bf p},{\bf q})=& \left\{g(\bx) \in  L_2(\cX): g(\bx)=\mu+\sum_{j \in \cJ} g_j(x_j) \geq 0\},\; |\cJ| \leq d_0,\right. \\
& \left.g_j \in B^{s_j}_{p_j,q_j}[0,1],\; j \in \cJ,\;\mE g_j(X_j)=0\right\},
\end{split}
\end{equation}
where $1 \leq p_j,q_j \leq \infty,\; (\frac{1}{p_j}-\frac{1}{2})_+ < s _j <r$ for all $j \in \cJ$, and the corresponding set of SpAM classifiers
$$
\cC_\cB(d_0,{\bf s}, {\bf p}, {\bf q})=\left\{\eta(\bx)=I\{g(\bx) \geq 0\}: g(\bx) \in \cG_{\cB}(d_0,{\bf s}, {\bf p}, {\bf q})\right\}.
$$

\begin{theorem} \label{th:besovupper}
Consider the SpAM  (\ref{eq:logistic})-(\ref{eq:SpAM}), where $g \in \cG_{\cB}(d_0,{\bf s}, {\bf p}, {\bf q})$, and assume $WRE(d_0,\bm,c_0)$-condition \ref{as:design} with $m_j=n^{\frac{1}{2s_j+1}},\;j=1,\ldots,d_0$ and $c_0$ that can be derived from the proof.

\begin{enumerate}
\item  Apply the sparse group Lasso classifier (\ref{eq:SpAMlasso}) with $\lambda=C_1 \sqrt{\frac{\ln d}{n}}$
and $\kappa=C_2 \sqrt{\frac{\ln n}{n}}$ for some $C_1, C_2>0$ that can be derived from the proof.
Then, 
$$
\sup_{\eta^* \in \cC_{\cB}(d_0, {\bf s}, {\bf p}, {\bf q})}\cE(\widehat{\eta}_{sgL},\eta^*) 
\le C \sqrt{\frac{1}{\nu_{sg}(d_0)}~\left(\frac{d_0 \ln d}{n}+
\ln n \sum_{j \in \cJ} n^{-\frac{2s_j}{2s_j+1}}\right)}
$$
for some $C>0$.

\item Apply the sparse group Slope classifier (\ref{eq:SpAMslope}) with $\lambda_j=C_1 \sqrt{\frac{\ln(de/j)}{n}},\;j=1,\ldots,d$
and $\kappa_\ell=C_2 \sqrt{\frac{\ln(ne/\ell)}{n}},\;\ell=1,\ldots,n$. 
Then, 
$$
\sup_{\eta^* \in \cC_{\cB}(d_0, {\bf s},{\bf p}, {\bf q})}\cE(\widehat{\eta}_{sgS},\eta^*) 
\le C \sqrt{\frac{1}{\nu_{sg}(d_0)}\left(\frac{d_0 \ln\left(\frac{de}{d_0}\right)}{n}+
\ln n \sum_{j \in \cJ} n^{-\frac{2s_j}{2s_j+1}}\right)}.
$$
\end{enumerate}
\end{theorem}

\begin{theorem} \label{th:besovlower}
Consider the SpAM  (\ref{eq:logistic})-(\ref{eq:SpAM}), where $d_0\ln(\frac{de}{d_0}) \leq n$ and $g \in \cG_{\cB}(d_0,{\bf s})$. Then,
$$
\inf_{\widetilde{\eta}}\sup_{\eta^* \in \cC_{\cB}(d_0, {\bf s},{\bf p}, {\bf q}),\mP_X}\cE(\widetilde{\eta},\eta^*) \geq \tilde{C} \sqrt{\frac{d_0 \ln\left(\frac{de}{d_0}\right)}{n}+
\sum_{j \in \cJ} n^{-\frac{2s_j}{2s_j+1}}}
$$
for some $\tilde{C}>0$.
\end{theorem}

Thus, sparse group Lasso and Slope wavelet-based SpAM classifiers with the tuning parameters from Theorem \ref{th:besovupper} are simultaneously nearly-optimal (up to log-factors) across the entire range of Besov classes including smooth and non-smooth functions. 
\ignore{
Moreover, similar to Sobolev classes, for sparse and moderately dense setups, $\widehat{\eta}_{sgS}$ is of the exact optimal order.
}

\section{Examples} \label{sec:example}
In this section we present experimental results for the developed SpAM classifiers applied to both simulated and real data.
Solving numerically (\ref{eq:sparsegroupLasso}) and (\ref{eq:sparsegroupSlope}) implies convex programming and various methods can be used. We could not find an available software for fitting logistic sparse group Slope and therefore tested  Lasso-type classifiers. In particular, we adapted the R package {\bf sparsegl} from CRAN based on the group-wise majorization-minimization algorithm
(see \cite{lchpm22} for details)\footnote{Note that the tuning parameters $\lambda$ and $\kappa$ of sprase group Lasso are defined slightly different in {\bf sparsegl.}}. Vanilla Lasso and group Lasso  were fitted as particular cases 
corresponding respectively to $\lambda=0$ and $\kappa=0$ in (\ref{eq:sparsegroupLasso}).   

\subsection{Simulations} \label{subsec:simulations}
We generated a sample of $n$  feature vectors $\bX_i=(X_{i1},\ldots,X_{id})^T \in \mathbb{R}^d$ according to
$$
X_{ij}=\frac{W_{ij}+U_i}{2},\;\;\;j=1,\ldots,d;\;\;i=1,\ldots, n,
$$
where $W_{ij}$ and $U_i$ are i.i.d. $U[0,1]$. Such a design results in correlations 0.5 between all features. 
We considered a SpAM with a logit function
$$
g(\bX)=g_1(X_1)+g_2(X_2)+g_3(X_3),
$$
where
$$
g_1(x_1)=|x_1-0.6|,\;\;\;g_2(x_2)=(2x_2-1)^2,\;\;\;g_3(x_3)=\frac{\sin(2\pi x_3)}{2-\sin(2\pi x_3)}
$$
   while other additive components are identically zero. The resulting vectors $g_1(\bX_1), g_2(\bX_2),g_3(\bX_3) \in \mathbb{R}^n$ were centred and scaled to have unit Euclidean norm. For a generated $\bX_i$, the class label $Y_i$ was sampled from a Bernoulli distribution $B(1,p(\bX_i)),\;i=1,\ldots,n$ with
$p(\bX_i)=\frac{e^{g(\bX_i)}}{1+e^{g(\bX_i)}}$.
We tried all combinations of $d=10, 30, 100$ and $n=30, 100, 300, 1000$ including, in particular, cases with $d \gtrsim n$.

We used cosine series basis, and applied Lasso, group Lasso and sparse group Lasso  (\ref{eq:sparsegroupLasso}) to find $\widehat{g}_{sgL}$ in (\ref{eq:g_sgL}). The tuning parameters were chosen by 10-fold cross-validation. The performance of the resulting SpAM classifiers (\ref{eq:SpAMlasso}) was measured by misclassification excess risk to compare the misclassification errors of the considered classifiers w.r.t. those of the Bayes (oracle) classifier $\eta^*(\bX)=I\{g(\bX) \geq 0\}$ on the independently generated test sets of size 100. The procedure was replicated 10 times for each combination of $d$ and $n$.

Table \ref{tab:simulations1} presents the misclassification excess risks on the test sets, the numbers of selected features (nonzero rows of the regression coefficient matrix $B$) and the overall numbers of selected coefficients (nonzero entries in $B$) averaged over 10 replications. 

\begin{table} 
    \caption{Average misclassification excess risks for various Lasso-type SpAM classifiers. The average numbers of selected features and of nonzero coefficients are given in parentheses.}
    \label{tab:simulations1}
    \begin{tabular}{cccc} 
    \hline
    $n$ & group Lasso & Lasso  & sparse group Lasso \\
    \hline
    \multicolumn{4}{c}{d=10}                           \\
    30   &  0.223       &  0.266        & 0.209        \\
         &  (3.7; 111)  &  (2.4; 4.3)   & (3.6; 75.4)  \\
    100  &  0.213       &  0.200        & 0.187        \\
         &  (7.6; 760)  &  (6.2; 16)    & (6; 213.4)   \\
    300  &  0.221       &  0.086        & 0.085        \\
         &  (8.9; 2670) &  (9.7; 45.1)  & (9.1; 108.9) \\
   1000  &  0.154       &  0.079        & 0.049        \\
         &  (9.8; 9800) &  (9.8; 49.3)  & (9.3; 94.6)  \\
                                                       \\        
    \multicolumn{4}{c}{d=30}                           \\
    30   &  0.258       &  0.319        & 0.251        \\
         &  (6.7; 201)  &  (1.9; 2.3)   & (4.6; 96.3)  \\
    100  &  0.232       &  0.178        & 0.177        \\
         &  (10.4; 1040)&  (14.9; 26)   & (10.6;3 3.9) \\
    300  &  0.192       &  0.067        & 0.070        \\
         &  (17.2; 5160)&  (16.3; 33.1) & (14.5; 40.6) \\
   1000  &  0.154       &  0.032        & 0.028        \\
         &  (9.8; 9800) &  (24.2; 77.6) & (15.6; 106.2)\\    
                                                       \\                                             
    \multicolumn{4}{c}{d=100}                          \\
    30   &  0.244       &  0.278        & 0.280        \\
         & (6; 180)     &  (9; 10.4)    & (4.9; 73.1)  \\
    100  &  0.257       &  0.207        & 0.192        \\
         & (19.7; 1370) &  (6.5; 7.5)   & (5.4; 7.7)   \\  
    300  &  0.203       &  0.101        & 0.092        \\
         & (25.8; 7740) &  (37.5; 53.2) & (38; 90.5)   \\ 
   1000  &  0.174       &  0.052        & 0.043        \\
         & (43.1; 43100)& (62.2; 142.0) & (37.8; 127.8)\\                 
    \hline
    \end{tabular}
    \end{table} 
    
Table \ref{tab:simulations1} demonstrates that a more flexible sparse group Lasso consistently outperforms vanilla Lasso and group Lasso.  
Misclassification excess risks for all methods converge to zero as $n$ increases for all $d$. 
%The convergence is faster for smaller $d$ since there are less %parameters to estimate.
Group Lasso converges more slowly and results in a larger number of selected coefficients, often leading to zero training errors. This behavior can be explained by the method keeping all coefficients of nonzero rows. 

It should be mentioned that while Lasso-type methods perform well for prediction, they may not be as successful in feature selection (support recovery) for strongly correlated features which are common for high-dimensional data (see, e.g., \cite{bg11}, Section 7.2). Thus, although in the vast majority of cases, all the procedures correctly selected $X_1$, $X_2$ and $X_3$, they often also included spurious features and coefficients especially for large $d$. Sparse group Lasso typically selected less features.

\subsection{Real-data example} \label{subsec:real_data}
To illustrate the performance of various Lasso-type classifiers on real data we considered the well-known benchmark email spam example from Example 1 of \cite{htf09}. The data consists of 4601 real and spam email messages, where for each email the outcome ({\bf mail} or {\bf spam}) and $d=57$ numeric attributes including relative frequencies of the most commonly occurring words and other email characters are available. The goal is to design a classifier for automatic spam detection based on this data. 

Following \cite{rllw07} we used only $n=300$ randomly selected samples for the training set to study the performance for relatively small $n$ (w.r.t. $d$) case. We used cosine series basis, and applied Lasso, group Lasso and sparse group Lasso  (\ref{eq:sparsegroupLasso}) to find $\widehat{g}_{sgL}$ in (\ref{eq:g_sgL}). The features were scaled first to have unit Euclidean norms and the test set was used as the hold-out set for choosing the tuning parameters. 
The performance of the resulting SpAM classifiers (\ref{eq:SpAMlasso}) was measured by misclassification errors on the remaining 4301 observations from the test set. For comparison, we
also present the results reported in \cite{rllw07} for their SpAM classifier, which is based on a certain additive smoothing procedure.

Table \ref{tab:results} provides misclassification errors for the test test, the subsets of selected features and the overall numbers of selected coefficients for the considered classifiers. Lasso-type classifiers resulted in smaller misclassification errors than SpAM of \cite{rllw07}. One realizes that sparse group Lasso classifier yielded the smallest misclassification error and both group-based classifiers outperformed vanilla Lasso. The latter disregarded the global row-wise sparsity of  the matrix $B$ and identified 32 active features (nonzero rows), although it resulted in the smallest overall number of nonzero coefficients.
Grouping coefficients from the same features allowed one to reduce the number of active features to 22, where both sparse group Lasso and group Lasso selected the same sets of features. While all entries of nonzero rows for group Lasso are nonzeroes (hence, $22 \times 300 = 6600$), sparse group Lasso, in addition, thresholded part of them as well.

    \begin{table}[h] 
    \caption{Misclassification errors and feature selection for various SpAM classifiers.}
    \label{tab:results}
    \begin{tabular}{lccc} 
    \hline
    Method & Error & \# Selected features & \# Nonzero coeff. \\
    \hline
    Lasso                                            &  0.0879 & 32 & 1703 \\
   \ignore{($\widehat{\kappa}=4\times10^{-5}$)}    &  & 
   $\{1-3,5-13,16-19,21,24,25,27,$    &             \\
   &  & $33,37,42,44-46,50,52,53,55-57\}$      & \\  
   group Lasso                                   &  0.0844 & 22  & 6600 \\
   \ignore{($\widehat{\lambda}=0.0029$)}  &  & $\{3,5,7,12,13,16,17,19,21,24,25,$ & \\
   & & $27,37,44-46,50,52-55,57\}$ & \\
   sparse group Lasso             & 0.0842  & 22  & 6389 \\
   \ignore{($\widehat{\lambda}=0.0030, \widehat{\kappa}=9\times10^{-6}$)}  & & 
   $\{3,5,7,12,13,16,17,19,21,24,25,$ & \\
   & & $27,37,44-46,50,52-55,57\}$ & \\
   SpAM  (\cite{rllw07})          & 0.1083  & 24  &   \\
   & & $\{4,6-10,14-22, 26,27,38, 53-58\}$ & \\
    \hline
    \end{tabular}
    \end{table} 

\medskip    
Summarizing, the results of this section indicate that Lasso-type classifiers demonstrate good performance across various types of data sets, and confirm the established theoretical findings. A more flexible sparse group Lasso overall outperforms vanilla Lasso which overlooks a grouped structure of coefficients, and group Lasso,
which fails to account for within-row sparsity and retains the entire nonzero rows. 

\section*{Acknowledgments}
The work was supported by the Israel Science Foundation (ISF), Grant ISF-1095/22. The author is grateful to Tomer Levy for valuable remarks. 

\appendix

\section{Appendix} \label{sec:appendix}
Throughout the proofs we use various generic positive constants, not necessarily the same each
time they are used even within a single equation.

\subsection{Proofs of the upper bounds} \label{subsec:appendix_upper}
%-- Theorems \ref{th:sgs_low_noise_classifier} and   \ref{th:nuclear_upper}} \label{sec:appendix_upper}

We start from a proof for a general upper bound and then apply it to correspondingly Sobolev, analytic and Besov classes of functions.

Consider a general logistic SpAM model (\ref{eq:model})-(\ref{eq:SpAM}), where $g_j \in \cF_j[0,1]$, and let $\{\psi_\ell\}_{l=0}^\infty$ be an orthonormal basis in $\cF_j$. Then,
$$
g(\bx)=\mu+\sum_{j=1}^d \sum_{\ell=1}^\infty \beta_{j\ell} \psi_\ell(x_j).
$$
Recall that $g_n$ is a truncated version of $g$:
$$
g_n(\bx)=\mu+\sum_{j=1}^d \sum_{\ell=1}^n \beta_{j\ell} \psi_\ell(x_j)
$$
and consider the corresponding truncated Bayes classifier $\eta^*_n=I\{g_n(\bx) \geq 0\}$.

Let $\widehat{g}$ and $\widehat{\eta}=I\{\widehat{g} \geq 0\}$ be respectively the logistic sparse group Lasso or Slope estimator of $g$ and the corresponding classifier. The misclassification excess risk $\cE(\widehat{\eta},\eta^*)$ can be decomposed as 
\be \label{eq:decomp}
\cE(\widehat{\eta},\eta^*)=\mE R(\widehat{\eta})-R(\eta^*_n)+R(\eta^*_n)-R(\eta^*)=\cE(\widehat{\eta},\eta^*_n)+R(\eta^*_n)-R(\eta^*).
\ee

\ignore{
The Kullback-Leibler divergence between two Bernoulli distributions with the odds $g_1$ and $g_2$ is 
$$
KL(g_1,g_2)=p_1\ln\frac{p_1}{p_2}+(1-p_1)\ln\frac{1-p_1}{1-p_2}=
\frac{e^{g_1}}{1+e^{g_1}}(g_1-g_2)+\ln(1+e^{g_2})-\ln(1+e^{g_1})
$$
Hence, the Kullback-Leibler divergence between the data distribution with the true $g_n$ and the empirical distribution generated by $\widehat{g}$ is
$d_{KL}(g_n,\widehat{g})=\int KL(g_n(\bx),\widehat{g}(\bx))d\mP_X$ and the corresponding Kullback-Leibler risk is $\mE d_{KL}(g,\widehat{g})$. From the well-known relations between the misclassification excess risk and the Kullback-Leibler risk we have
$$
\cE(\widehat{\eta},\eta^*_n) \leq \sqrt{2 \mE d_{KL}(g_n,\widehat{g})}
$$
(\cite{y99}). 
}

The logistic loss $\ell(y, g_n(\bx))$ in (\ref{eq:logloss}) is Lipschitz. In addition, one can show that for the bounded support $\cX$, $\mE\left(\ell(Y,g_{1n}(\bX)-\ell(Y,g_{2n}(\bX))\right) \geq C ||g_{1n}(\bx)-g_{2n}(\bx)||_{L_2(\mP_X)}^2$ for some $C>0$ and any pair $(g_{1n}, g_{2n})$. Furthermore, assume that it is possible to approximate $d_0$ ``approximately sparse'' nonzero rows $B_{j\cdot},\;j \in \cJ$ of the regression coefficients matrix $B$ by truly sparse vectors $B'_{j\cdot}$ 
with the properly chosen numbers of nonzero entries
$m_j$'s such that $\kappa |B_{j\cdot}-B'_{j\cdot}|_1=\kappa \sum_{\ell=1}^n |\beta_{j\ell}-\beta'_{j\ell}|
\leq \ln n~ m_j n^{-1}$ for Lasso and $\sum_{\ell=1}^n \kappa_j |B-B'|_{j(l)} \leq  \ln(ne/m_j) m_j n^{-1}$ for Slope for $\kappa$'s given in the theorems. Similar to $B$, complete the matrix $B'$ with zero rows for $j \notin \cJ$. One can then exploit the results of \cite{la23} (Section 3.2 and Remark 8) for sparse group Lasso and Slope to get under the $WRE(d_0,\bm,c_0)$-condition \ref{as:design}\footnote{Although the original assumption required in \cite{la23} is somewhat stronger, it can be relaxed to the WRE condition \ref{as:design}.}
\be \label{eq:general_lasso}
\cE(\widehat{\eta}_{sgL},\eta^*_n) 
\le C \sqrt{\frac{1}{\nu_{sg}(d_0)}~\frac{d_0 \ln d+ \ln n \sum_{j \in \cJ} m_j }{n}}
\ee
and
\be \label{eq:general_slope}
\cE(\widehat{\eta}_{sgS},\eta^*_n) 
\le C \sqrt{\frac{1}{\nu_{sg}(d_0)}~\frac{d_0 \ln\left(\frac{ de}{d}\right)+ \sum_{j \in \cJ} m_j \ln\left(\frac{ne}{m_j}\right)}{n}}.
\ee

To bound the second term in the RHS of (\ref{eq:decomp}) note that
\be \label{eq:tails}
R(\eta^*_n)-R(\eta^*) \leq 
%\sqrt{2 KL(g,g_n)} \leq 
2||g_n-g||_{L_2(\mP_X)} \leq C \sqrt{\sum_{j \in \cJ} \sum_{\ell \geq n} |\beta_{j\ell}|^2}
\ee
(\cite{z04a}). We will show that this term is negligible w.r.t.  $\cE(\widehat{\eta},\eta^*_n)$.

We now adapt the presented general scheme for Sobolev, analytic and Besov classes by finding the corresponding matrices $B'$.

\subsubsection{Proof of Theorem \ref{th:sobolevupper}} \label{subsubsec:sobolevupper}
For $g_j \in \widetilde{H}^{s_j}[0,1],\;s_j>\frac{1}{2}$ construct the corresponding sparse approximations $B'_{j\cdot}$ for nonzero rows  $B_{j\cdot}$ of $B$ as follows. 
Take $m_j=n^{\frac{1}{2s_j+1}}$ and set $\beta'_{j\ell}=\beta_{j\ell}$ for the first $m_j$ entries of $B_{j\cdot}$ and zeroes for others. Note that for $g_j \in \widetilde{H}^{s_j}[0,1], s_j>1/2$
the cosine series coefficients satisfy $|\beta|_{j\ell} \leq C \ell^{-(s_j+\frac{1}{2})}$. 
For Lasso with $\kappa=C_2 \sqrt{\frac{\ln n}{n}}$ 
%with $C_2$ given in the proofs of Levy and Abramovich (2023)
we then have
$$
\kappa |B_{j\cdot}-B'_{j\cdot}|_1 = \kappa \sum_{\ell=m_j+1}^n |\beta|_{j\ell} \leq  C \kappa~ m_j^ {-s_j+\frac{1}{2}} \leq C \sqrt{\ln n}~ n^{-\frac{2s_j}{2s_j+1}} \leq C \frac{\ln n}{n} m_j.
$$
Similarly, for Slope with $\kappa_j=C_2 \sqrt{\frac{\ln(ne/j)}{n}}$,
\begin{equation} \nonumber
\begin{split}
\sum_{\ell=1}^n \kappa_\ell |B-B'|_{j(\ell)}& =\sum_{\ell=m_j+1}^n
\kappa_l |\beta|_{j\ell} \leq C \kappa_{m_j} ~ m_j^{-s_j+\frac{1}{2}} \leq C \sqrt{\ln\left(\frac{ne}{m_j}\right)}n^{-\frac{2s_j}{2s_j+1}}\\ 
& \leq C \frac{\ln\left(\frac{ne}{m_j}\right)}{n}m_j.
\end{split}
\end{equation}

Hence, applying the results of \cite{la23} (Section 3.2 and Remark 8)
with $\lambda$'s and $\kappa$'s in Theorem \ref{th:sobolevupper} and constants $C_1, C_2$ given in their proofs, (\ref{eq:general_lasso}) implies
\begin{equation} \nonumber
\begin{split}
\cE(\widehat{\eta}_{sgL},\eta^*_n) 
& \le C \sqrt{\frac{1}{\nu_{sg}(d_0)}~\frac{d_0 \ln d+ \ln n \sum_{j \in \cJ} m_j }{n}} \\
& \leq
C \sqrt{\frac{1}{\nu_{sg}(d_0)}~\left(\frac{d_0 \ln d}{n}+
\ln n \sum_{j \in \cJ} n^{-\frac{2s_j}{2s_j+1}}\right)}.
\end{split}
\end{equation}
Similarly, by (\ref{eq:general_slope})
\be \nonumber 
\begin{split}
\cE(\widehat{\eta}_{sgS},\eta^*_n) 
& \le C \sqrt{\frac{1}{\nu_{sg}(d_0)}~\frac{d_0 \ln\left(\frac{ de}{d_0}\right)+ \sum_{j \in \cJ} m_j \ln\left(\frac{ne}{m_j}\right)}{n}} \\
&\leq
C \sqrt{\frac{1}{\nu_{sg}(d_0)}~\left(\frac{d_0\ln\left(\frac{ de}{d_0}\right) }{n}+
\ln n \sum_{j \in \cJ} n^{-\frac{2s_j}{2s_j+1}}\right)}.
\end{split}
\ee

Finally, for the tail sums in (\ref{eq:tails}) we have
$$
\sqrt{\sum_{j \in \cJ}\sum_{\ell > n} \beta_{j\ell}^2} \leq C \sqrt{\sum_{j \in \cJ} n^{-2s_j}}=o(\cE(\widehat{\eta},\eta^*_n)).
$$

\subsubsection{Proof of Theorem \ref{th:analyticupper}} \label{subsubsec:analyticupper}
For an analytic $g_j \in \cA^{\balpha_j}[0,1]$, the cosine coefficients $|\beta|_{j\ell} \leq C e^{-\alpha_j \ell}$. Truncate nonzero rows of the regression matrix $B$ at $m_j=\frac{1}{\alpha_j}\ln \frac{n}{\ln n}, \; j \in \cJ$. 
We have
$$
\kappa |B_{j\cdot}-B'_{j\cdot}|_1=\kappa \sum_{\ell=m_j+1}^n |\beta|_{j \ell} \leq C \kappa e^{-\alpha_j m_j} 
\leq C \sqrt{\frac{\ln n}{n}} \leq C \frac{\ln n}{n} m_j
$$
for Lasso and
$$
\sum_{\ell=1}^n \kappa_\ell |B-B'|_{j(\ell)}=\sum_{\ell=m_j+1}^n \kappa_\ell |\beta|_{j\ell} \leq C \kappa_{m_j} e^{-\alpha_j m_j} \leq C \frac{\ln\left(\frac{ne}{m_j}\right)}{n} m_j
$$
for Slope.

Going along the lines of the proof for Sobolev classes, (\ref{eq:general_lasso}) and (\ref{eq:general_slope}) yield then
\begin{equation} \nonumber
\begin{split}
\cE(\widehat{\eta}_{sgL},\eta^*_n) & 
\le C \sqrt{\frac{1}{\nu_{sg}(d_0)}~\frac{d_0 \ln d+ \ln n \sum_{j \in \cJ} m_j }{n}} \\
& \le C \sqrt{\frac{1}{\nu_{sg}(d_0)}~\left(\frac{d_0 \ln d}{n}+
\sum_{j \in \cJ}\frac{1}{\alpha_j}~ \frac{\ln n}{n}\right)}
\end{split}
\end{equation}
and
\be \nonumber
\begin{split}
\cE(\widehat{\eta}_{sgS},\eta^*_n) 
& \le C \sqrt{\frac{1}{\nu_{sg}(d_0)}~\frac{d_0 \ln\left(\frac{ de}{d_0}\right)+ \sum_{j \in \cJ} m_j \ln\left(\frac{ne}{m_j}\right)}{n}} \\
& \le C \sqrt{\frac{1}{\nu_{sg}(d_0)}~\left(\frac{d_0 \ln\left(\frac{de}{d_0}\right)}{n}+
\sum_{j \in \cJ}\frac{1}{\alpha_j}~ \frac{\ln n}{n}\right)}.
\end{split}
\ee
Evidently, $\sqrt{\sum_{j \in \cJ}\sum_{\ell > n} \beta_{j\ell}^2}=o\left(\cE(\widehat{\eta}_{sgS},\eta^*_n)\right)$.

\subsubsection{Proof of Theorem \ref{th:besovupper}} \label{subsubsec:besovupper}
Recall that for $g_j \in B^{s_j}_{p_j,q_j}[0,1]$ we consider orthonormal wavelet series $\{\psi_\ell\}_{l=0}^\infty$ of regularity $r>s_j$, where $\ell=2^h+k,\;h \geq 0,\;0 \leq k \leq 2^h-1$, and $|\beta|_{j(\ell)} \leq C \ell^{-(s_j+\frac{1}{2})}$ (see (\ref{eq:weaklp})). For a given nonzero row of $B$ design its sparse approximation $B'_{j\cdot}$ as follows: keep $m_j=n^{\frac{1}{2s_j+1}}$ coefficients of $B_{j\cdot}$ with the largest absolute values and zero others.

Hence,
$$
\kappa |B_{j\cdot}-B'_{j\cdot}|_1 = \kappa \sum_{\ell=m_j+1}^n|\beta|_{j(\ell)} \leq  C \kappa~ m_j^ {-s_j+\frac{1}{2}} \leq C \sqrt{\ln n}~ n^{-\frac{2s_j}{2s_j+1}} \leq C \frac{\ln n}{n} m_j
$$
for Lasso and
\begin{equation} \nonumber
\begin{split}
\sum_{\ell=1}^n \kappa_\ell |B-B'|_{j(\ell)}&=\sum_{\ell=m_j+1}^n
\kappa_l |\beta|_{j(\ell)} \leq C \kappa_{m_j} ~ m_j^{-s_j+\frac{1}{2}} \leq C \sqrt{\ln\left(\frac{ne}{m_j}\right)}~ n^{-\frac{2s_j}{2s_j+1}} \\
& \leq C \frac{\ln\left(\frac{ne}{m_j}\right)}{n} m_j
\end{split}
\end{equation}
for Slope. To complete the proof one repeats the arguments from the proof for Sobolev classes.

\subsection{Proofs of the lower bounds (Theorems \ref{th:sobolevlower}, \ref{th:analyticlower} and \ref{th:besovlower})} \label{subsec:appendix_lower}
To prove the lower bounds for the minimax misclassification excess risk we utilize the general approach developed in \cite{y99}.

Assume that the unknown logit function $g(\bx)$ belongs to some class of functions $\cG$ with a metric entropy 
$\cm(\cG,\epsilon)$ -- the logarithm of the $\epsilon$-packing number of $\cG$ w.r.t. $L_2(\mP_X)$-norm. Consider the corresponding class $\cP$ for 
$p(\bx)$: $\cP=\{p: p=\frac{e^g}{1+e^g},\; g \in \cG\}$.
Let $\epsilon^*$ satisfy 
\be \label{eq:entropy1}
n\epsilon^{*2} \gtrsim \cm(\cP,\epsilon^*).
\ee
Then if the class $\cP$ is rich enough in terms of $\cm(\cP,\epsilon)$ to satisfy the assumptions given in \cite{y99}, the minimax misclassification excess risk
\be \label{eq:entropy-lower_bound}
\inf_{\widetilde{\eta}} \sup_{p \in \cP,\mP_X} \cE(\widetilde{\eta},\eta^*) \asymp \epsilon^*,
\ee
where recall that the Bayes classifier $\eta^*(\bx)=I\{p(\bx) \geq 1/2\}$.

By Taylor expansion, one immediately verifies that for any probability functions $p_1(\bx)$ and $p_2(\bx)$ and their logit functions $g_1(\bx), g_2(\bx)$,
$$
||p_2-p_1||_{L_2(\mP_X)} \leq \frac{1}{4} ||g_2-g_1||_{L_2(\mP_X)}.
$$
Hence, $\cm(\cP,\epsilon) \leq \cm(\cG, 4\epsilon)$.

For the considered general SpAM setup with $d_0 < d/4$ by Lemma 4 of Raskutti {\em et al.} (2012),
\be \label{eq:entropy2}
\cm(\cG,4\epsilon) \leq C \left(d_0 \ln \frac{de}{d_0} + \sum_{j \in \cJ}\cm\left(\cG_j,\frac{4\epsilon}{\sqrt{d_0}}\right)\right),
\ee
where for the considered Sobolev, analytic and Besov classes $\cG_j$'s we have $\cm(\widetilde{H}^{s_j}[0,1],\epsilon) \simeq
\cm(B^{s_j}_{p_j,q_j}[0,1],\epsilon)  
\simeq
\left(\frac{1}{\epsilon}\right)^{\frac{1}{s_j}}$ and 
$\cm(\cA^{\alpha_j}[0,1],\epsilon) \simeq \frac{1}{\alpha_j} \ln \left(\frac{1}{\epsilon}\right)$ (e.g., \cite{y99, rwy12}).
Moreover, \cite{y99} showed that the required conditions on metric entropies to apply (\ref{eq:entropy-lower_bound}) are satisfied for these classes. 

Solving (\ref{eq:entropy1}) and (\ref{eq:entropy2}) after some algebra implies that the suitable solutions are $\epsilon^{*2} =C \left(\frac{d_0\ln \frac{de}{d_0}}{n} + \sum_{j \in \cJ}n^{-\frac{2s_j}{2s_j+1}}\right)$ and $\epsilon^{*2} =
C\left(\frac{d_0\ln \frac{de}{d_0}}{n} + \sum_{j \in \cJ} \frac{1}{\alpha_j} \frac{\ln n}{n}\right)$ respectively.

For the dense case $d/4 \leq d_0 \leq d$ (hence, $d_0 \sim d$), 
the model selection error is negligible w.r.t. nonparametric estimation term.
Evidently, no classifier can perform better than an oracle who knows 
the subset of truly active features $\cJ$. Therefore, for Sobolev and Besov classes in this case
$d_0 \ln\frac{de}{d_0} \sim d_0 
\lesssim \sum_{j \in \cJ} n^{\frac{1}{2s_j+1}}$ and $\epsilon^{*2}= C \sum_{j \in \cJ} n^{-\frac{2s_j}{2s_j+1}}$. Similarly, for analytic functions, $d_0 \ln\frac{de}{d_0} \lesssim \ln n  \sum_{j \in \cJ}\frac{1}{\alpha_j}$ and
$\epsilon^{*2}=C\frac{\ln n}{n}  \sum_{j \in \cJ}\frac{1}{\alpha_j}$.

Applying (\ref{eq:entropy-lower_bound}) completes the proof.

\end{document}